\documentclass[12pt]{article}
\font\bigbf=cmbx12 at 16pt

\usepackage[english,activeacute]{babel}

\usepackage{amsfonts,amssymb,amstext}
\usepackage{mathrsfs}
\usepackage{fancyhdr}
\usepackage{amssymb, bbm}                          
\usepackage{amssymb,amsmath,babel,graphicx,epsfig} 
\usepackage{color}

\def\R {{\rm I}\hskip -0.85 mm{\rm R }}

\def\N {{\rm I}\hskip -0.85 mm{\rm N }}

\def\n{{\bf n}}

\def\u{{\bf u}}

\def\w{{\bf w}}

\def\wb4b{{\bf W}_b^{1,4}(\Om)}

\def\Om{\Omega}

\def \beq {\begin{equation}}
\def \eeq {\end{equation}}
\def \ba {\begin{array}}
\def \ea {\end{array}}
\def \dis {\displaystyle}
\def \ecart{\noalign{\medskip}}

\def\espan{\mathop{\rm span}\nolimits}

\def\u3d{{\bf U}}
\def\w3d{{\bf W}}

\def\n3d{{\bf N}}

\def\csu{{\cal S}_1}\def\csk{{\cal S}_k}

\def\ug{u(\gamma) }
\def\ugw{u_W(\gamma)}
\def\ugz{u_Z(\gamma)}
\def\ab{\overline{a}}

\def\w4{W_b^{1,4}(\Omega)}

\def\l2{L_0^2(\Omega,\partial_3)}

\def\R{\mathbb{R}}

\newtheorem{theorem}{Theorem}[section]
\newtheorem{lemma}[theorem]{Lemma}
\newtheorem{proposition}[theorem]{Proposition}
\newtheorem{corollary}[theorem]{Corollary}

\newtheorem{remark}[theorem]{Remark}

\def\blackbox{\leavevmode\vrule height 5pt width 4pt depth 0pt\relax}
\newenvironment{proof}{\begin{trivlist}
                       \item[]\hspace{0cm}{\bf Proof:}\hskip -5pt
                       \hspace{0cm} }{\hfill $\blackbox$
                     \end{trivlist}}

\title{\bigbf An intrinsic Proper Generalized Decomposition 
\\
for parametric symmetric elliptic problems}

\begin{document}

\author{
 M. Aza\"\i ez \thanks{I2M, IPB (UMR CNRS 5295), Universit\'e de Bordeaux, 33607 Pessac (France).}
\qquad 
  F. Ben Belgacem \thanks{Universit\'e Paris Sorbonne, UTC, EA 2222, Laboratoire de Math\'ematiques Appliqu\'ees de Compi\`egne,  60205 Compi\`egne (France).}
 \qquad
  J. Casado-D\'{\i}az \thanks{Departamento EDAN \& IMUS, Universidad de Sevilla, C/Tarfia, s/n, 41012 Sevilla (Spain).}
\\
T. Chac\'on Rebollo \thanks{I2M, IPB (UMR CNRS 5295), Universit\'e de Bordeaux, 33607 Pessac (France), and Departamento EDAN \& IMUS, Universidad de Sevilla, C/Tarfia, s/n, 41012 Sevilla (Spain).}
\qquad 
  F. Murat \thanks{Laboratoire Jacques-Louis Lions, Bo\^{\i}te courrier 187,
Universit\'e Pierre et Marie Curie (Paris VI),
75252 Paris cedex 05 (France). }
}
\maketitle

\abstract{We introduce in this paper a technique for the reduced order approximation of parametric symmetric elliptic partial differential equations. 
For any given dimension, we prove the existence of an optimal subspace of at most that dimension which realizes 
the best approximation in mean of the error with respect to the parameter
in the quadratic norm associated to the elliptic operator
between the exact solution and the Galerkin solution calculated on the subspace. This is analogous to the best approximation property of the Proper Orthogonal Decomposition (POD) subspaces, excepting that in our case the norm is parameter-depending, and then  the POD optimal sub-spaces cannot be characterized by means of a spectral problem. We apply a deflation technique to build a series of approximating solutions on finite-dimensional optimal subspaces, directly in the on-line step. We prove that the partial sums converge to the continuous solutions 
in mean quadratic elliptic norm.}

\clearpage

\section{Introduction}
The Karhunen-Lo\`eve's expansion  (KLE) is a widely used tool, that provides a reliable procedure for a low dimensional representation of spatiotemporal signals  (see \cite{Ghanem, loeve}).  It is referred to as the  principal components analysis (PCA) in statistics (see \cite{Hotelling, PCA, Pearson}), or called singular value decomposition (SVD) in linear algebra (see \cite{golub}). It is named  the proper orthogonal decomposition (POD) in mechanical computation, where it is also widely used (see \cite{Berkoz}). Its use allows large savings of computational costs, and make affordable the solution of problems that need a large amount of solutions of parameter-depending Partial Differential Equations (see \cite{aeroelastic,  cordier, lumley, lass, Pearson, modelredb, modelred,  volkwein}). 
\par
However the computation of the POD expansion requires to know the function to be expanded, or at least its values at the nodes of a fine enough net. This makes it rather expensive to solve parametric elliptic Partial Differential Equations (PDEs), as it requires the previous solution of the PDE for a large enough number of values of the parameter (\lq\lq snapshots\rq\rq) (see \cite{kavo}), even if these can be located at optimal positions (see \cite{optimalsnap}).  Galerkin-POD strategies are well suited to solve parabolic problems, where the POD basis is obtained from the previous solution of the underlying elliptic operator (see \cite{kunish, muller}).

An alternative approach is the Proper Generalized Decomposition that iteratively computes a tensorized representation of the parameterized PDE, that separates the parameter and the independent variables, introduced in \cite{ammar}.  It has been interpreted as a Power type Generalized Spectral Decomposition (see \cite{nouy1, nouy2}). It has experienced a fast development, being applied to the low-dimensional tensorized solution of many applied problems. The mathematical analysis of the PGD has experienced a relevant development in the last years. The convergence of a version of the PGD for symmetric elliptic PDEs via minimization of the associated energy has been proved in \cite{BrisLeMad}. Also,  in \cite{falnou} the convergence of a recursive approximation of the solution of a linear elliptic PDE is proved, based on the existence of optimal subspaces of rank 1 that minimize the elliptic norm of the current residual. 

The present paper is aimed at the direct determination of a variety of reduced dimension for the solution of parameterized symmetric elliptic PDEs. We intend to on-line determine an optimal subspace of given dimension that yields the best approximation in mean (with respect to the parameter) of the error (in the quadratic norm associated to the elliptic operator) between the exact solution and the Galerkin solution calculated on the subspace.  The optimal POD sub-spaced can no longer be characterized by means of a spectral problem for a compact self-adjoint operator (the standard POD operator) and thus the spectral theory for compact self-adjoint operators does no apply.  We build recursive approximations on finite-dimensional optimal subspaces by minimizing the mean quadratic error of the current residual, similar to the one introduced in \cite{falnou}, that we prove to be strongly convergent in the \lq\lq natural\rq\rq mean quadratic elliptic norm. The method shares some properties of PGD and POD expansions: It builds a tensorized representation of the parameterized solutions, by means of optimal subspaces that minimize the residual in mean quadratic norm.
\par
The paper is structured as follows: In Section 2 we state the general problem of finding optimal subspaces of a given dimension. We prove in Section 3 that there exists a solution for 1D optimal subspaces, characterized as a maximization problem with a non-linear normalization restriction. We extend this existence result in Section 4 to general dimensions. Finally, in Section 5 we use the results in Sections 3 and 4 to
build a deflation algorithm to approximate the solution of a parametric family of elliptic problems and we show the convergence.
  \clearpage 

\section{Statement of the problem}

Let $H$ be a separable Hilbert space endowed with the scalar product $ (\cdot,\cdot)$.
The related norm is denoted by $\|\cdot \|$.\par
 We  denote by $B_s(H)$ the space of bilinear, symmetric and continuous forms in $H$.\par
Assume given a measure space  $(\Gamma, {\cal B},\mu)$, with standard notation, so that $\mu$
is $\sigma$-finite. \par
Let $a\in L^\infty(\Gamma,B_s(H);d\mu)$ be
such that there exists $\alpha>0$ satisfying
\begin{equation}\label{eq:propab}
\alpha\, \|u\|^2 \le a(u,u;\gamma),\quad \forall u \in H,\ d\mu\hbox{-a.e. }\gamma\in\Gamma.
\end{equation}
For  $\mu-$a.e $\gamma\in \Gamma$, the bilinear form $a(\cdot,\cdot;\gamma)$ determines a norm uniformly equivalent to the norm $\|\cdot \|$. Moreover,   $\ab\in B_s(H)$ defined by
\begin{equation} \label{defab}
\ab(v,w) = \int_\Gamma a(v(\gamma),w(\gamma);\gamma)\, d\mu(\gamma),\quad \forall  v,\, w \in L^2(\Gamma,H;d\mu)
\end{equation}
defines an inner product in $H$ which generates a norm  equivalent to the standard one in $L^2(\Gamma,H;d\mu)$.\par\medskip
Let be  given  a data function $f \in L^2(\Gamma, H';d\mu)$.  We are interested in the variational problem:
\begin{equation}\label{eq:pbg} \hbox{ Find }\ \ug \in H\ \hbox{ such that }\quad
 a(u(\gamma),v;\gamma)=\langle f(\gamma),v \rangle,\quad \forall v \in H,\ d\mu\hbox{-a.e. }\gamma\in \Gamma,
 \end{equation}
where $\langle \cdot, \cdot \rangle$ denotes the duality pairing between $H'$ and $ H$.\par
By Riesz representation theorem, problem  \eqref{eq:pbg} admits a unique solution for $d\mu$-a.e. $\gamma \in \Gamma$. On the other hand, we claim that $\tilde u$ solution of
\beq\label{eq:pbg2}\tilde u\in L^2(\Gamma,H;d\mu),\qquad \bar a(\tilde u, \bar v)=\int_{\Gamma}\langle f(\gamma),\bar v(\gamma)\rangle\, d\mu(\gamma),\quad \forall\, \bar v\in L^2(\Gamma,H;d\mu),\eeq
also satisfies (\ref{eq:pbg}): Indeed taking $\bar v=v\chi_B$, with $v\in H$ fixed and $B\in {\cal B}$ arbitrary, implies that there exists a subset $N_v\in {\cal B}$ with $\mu(N_v)=0$ such that
$$a(\tilde u(\gamma),v;\gamma)=\langle f(\gamma),v \rangle,\quad\forall\,\gamma\in \Gamma\setminus N_v.$$
The separability of $H$ implies that $N_v$ can be chosen independent of $v$, which proves the claim. By the uniqueness of the solution of (\ref{eq:pbg}) this shows that
\beq\label{dealu} \tilde u=u\quad d\mu\hbox{-a.e. }\gamma\in \Gamma.\eeq 
This proves that $u$ belongs to $L^2(\Gamma,H;d\mu)$ and provides an equivalent definition of $u$. Namely, that $u$ is  the solution of  (\ref{eq:pbg2}).\par
Given a closed subspace $Z $ of $H$, let us denote by $u_Z(\gamma)$ the solution of the Galerkin approximation of problem (\ref{eq:pbg}) on $Z$, which reads as 
 \begin{equation}\label{eq:pbgw}
 u_Z(\gamma) \in Z,\quad a(\ugz,z;\gamma)=\langle f(\gamma),z \rangle,\quad \forall z \in Z,\ d\mu\hbox{-a.e. }\gamma\in \Gamma,
\end{equation}
or equivalently as
\begin{equation}\label{eq:pbgw2} u_Z\in L^2(\Gamma,Z;d\mu),\qquad \bar a(u_Z,z)=\int_\Gamma\langle f(\gamma),z(\gamma)\rangle\, d\mu(\gamma),\quad \forall\, z\in L^2(\Gamma,Z;d\mu).
 \end{equation}
For every $k\in\N$, we intend to find the best subspace $W $ of $H$ of dimension smaller than or equal to $k$  that minimizes the mean error between $\ug$ and $\ugw$. That is, $W$ solves
\begin{equation}\label{eq:pb1d}
 \qquad \min_{Z \in \csk}  \bar a(u-u_Z,u-u_Z),
\end{equation}
where ${\cal S}_k$ is the family of subspaces of $H$ of dimension smaller than or equal to $k$. This problem will be  proved to have a solution in Sections 3 and 4. We will then use this result to approximate the solution $u$ of problem (\ref{eq:pbg}) by a deflation algorithm.

\par\medskip
To finish this section we  provide some equivalent formulations of the problem. First we 
observe that
\begin{proposition} \label{procaruz} For every closed subspace $Z\subset H$, the function $u_Z$ defined by (\ref{eq:pbgw2}) is also the unique solution of 
\beq\label{caruz}\min_{z\in L^2(\Gamma,Z;d\mu)}\bar a(u-z,u-z).\eeq
Moreover, for $d\mu$-a.e. $\gamma\in \Gamma$, the vector $u_Z(\gamma)$ is the solution of 
\beq\label{caruz2} \min_{z\in Z} a(u(\gamma)-z,u(\gamma)-z;\gamma).\eeq
\end{proposition}
\begin{proof} It  is a classical property of the Galerkin  approximation of the variational formulation of  linear elliptic problems that $u_Z$ satisfies  (\ref{caruz}). Indeed, the symmetry of $\bar a$ gives
$$\bar a(u-z,u-z)=\bar a(u-u_Z,u-u_Z)+2\bar a(u-u_Z,u_Z-z)+\bar a(u_Z-z,u_Z-z),$$
for every $z\in L^2(\Gamma,H;d\mu)$, where by (\ref{eq:pbg2}), (\ref{dealu}) and (\ref{eq:pbgw2}) the second term on the right-hand side vanishes, while the third one is nonnegative. This proves (\ref{caruz}).\par
The proof of (\ref{caruz2}) is the same by taking into account (\ref{eq:pbg}) and (\ref{eq:pbgw}) instead of (\ref{eq:pbg2}) and (\ref{eq:pbgw2}).
\end{proof}\par
As a consequence of Proposition \ref{procaruz} and definition (\ref{defab}) of $\bar a$, we have
\begin{corollary}  A space $W\in S_k$ is a solution of (\ref{eq:pb1d}) if and only if it is a solution of
\beq\label{eq:pb1db} \min_{Z\in S_k}\min_{z\in L^2(\Gamma,Z;d\mu)}\bar a(u-z,u-z).\eeq
Moreover
\beq\label{eq:pb1dbbis} \min_{Z\in S_k}\min_{z\in L^2(\Gamma,Z;d\mu)}\bar a(u-z,u-z)=
\min_{Z\in S_k}\int_\Gamma\min_{z\in Z} a(u(\gamma)-z,u(\gamma)-z;\gamma)d\mu(\gamma).\eeq

\end{corollary}

\begin{remark}
Optimization problem \eqref{eq:pb1db} is reminiscent of the Kolmogorov  $k$-width related to the best approximation of the manifold $(u(\gamma))_{\gamma\in \Gamma}$ by subspaces in $H$  with dimension $k$ as  presented in {\em \cite{Maday}}.
In  the present minimization problem, we use   the norm of $L^2(\Gamma,H,d\mu)$ instead of
the norm of $L^\infty(\Gamma,H,d\mu)$ as used there. The minimization problem in {\em \cite{Maday}}  can indeed be written as  
$$
\min_{Z \in \mathcal S_{ k}}\, \hbox{\rm ess}\hskip-6pt\sup_{\gamma \in \Gamma\quad}\min_{z\in Z} \,
a(\ug-z,\ug-z;\gamma),
$$
if one uses $a(\cdot,\cdot;\gamma)$ as the inner product in $H$.\\
\end{remark}

For a function $v\in L^2(\Gamma,V;d\mu)$, we denote by $R(v)$ the vectorial space spanned by $v(\gamma)$ when $\gamma$ belongs to $\Gamma$; more exactly, taking into account that $v$ is only defined up to sets of zero measure, the correct definition of $R(v)$ is given by
\begin{equation}\label{DefRv} R(v)=\bigcap_{\mu( N)=0}\hbox{Span }\big\{v(\gamma):\ \gamma\in\Gamma\setminus  N\big\}.\end{equation}
\par
Taking into account (\ref{eq:pb1db}), a new formulation of (\ref{eq:pb1d}) is given by
\begin{proposition}\label{pbeq1} If $W$ is a solution of (\ref{eq:pb1d}), then $u_W$ is a solution of
\beq\label{caruW} \min_{v\in L^2(\Gamma,H;d\mu)\atop
{\rm dim}\,R(v)\leq k}\bar a(u-v,u-v).\eeq
Reciprocally, if $\hat u$ is a solution of (\ref{caruW}), then $R(\hat u)$ is a solution of (\ref{eq:pb1d}) and $\hat u=u_{R(\hat u)}$.
\end{proposition}
\par
 The next proposition provides another formulation for (\ref{eq:pb1d}) which depends on $f$ and not on the solution $u$ of (\ref{eq:pbg}).
\begin{proposition} \label{pr:equiv} The subspace $W \in \csk$ solves problem \eqref{eq:pb1d}
 if and only if it is a solution  of the problem
\begin{equation}\label{eq:pb1dr}
\max_{Z \in \csk}  \int_\Gamma \langle f(\gamma),u_Z(\gamma)\rangle \, d\mu(\gamma).\end{equation}
\end{proposition}

\begin{proof} As in the proof of the first part of Proposition \ref{procaruz}, one deduces from
 (\ref{eq:pbg2}), (\ref{dealu}) and (\ref{eq:pbgw2}) that
$$
\bar a(u-u_Z,z)=0,\quad \forall z \in L^2(\Gamma,Z;d\mu).
$$
Using the symmetry of $\bar a$, we then have
$$\begin{array}{l}
\bar a(u-u_Z,u-u_Z) 
=\bar a(u,u)- a(u_Z,u)=\bar a(u,u)-\bar a(u_Z,u_Z)
\\ \noalign{\medskip}\displaystyle=  \bar a(u,u)- \int_\Gamma \langle f(\gamma),u_Z(\gamma)\rangle \, d\mu(\gamma).
\end{array}$$
Thus $W$ solves (\ref{eq:pb1d}) if and only if it solves (\ref{eq:pb1dr}).
\end{proof}
\
\begin{remark}
In  {\em \cite{falnou}} a problem similar to \eqref{eq:pb1d} has been studied, namely
\begin{equation}\label{eq:pbfalnou}
(P_k)' \qquad \min_{Z \in \csk}  \int_\Gamma (\ug-\ugz,\ug-\ugz)_H\,d\mu(\gamma),
\end{equation}
where $(\cdot,\cdot)_H$ is an inner product on $H$. In this case a solution of $(P_k)'$ is the space generated by the first $k$ eigenfunctions of the POD operator ${\cal P}: H \mapsto H$, which is given by
$$
{\cal P}(v) =\int_\Gamma (v, u(\gamma))_H \, u(\gamma)\, d \mu(\gamma),\quad \forall v \in H.
$$
\par In the present case, due the dependence of $a$ with respect to $\gamma$, it does not seem that the problem can be reduced to a spectral problem. As an example, we consider the case $k=1$. Then problem
(\ref{caruW}) can be written as
\beq\label{caruWbb} \min_{v\in H,\,\varphi\in L^2(\Gamma;d\mu)}
\int_\Gamma a(u(\gamma)-\varphi(\gamma) v,u(\gamma)-\varphi(\gamma) v;\gamma)d\mu(\gamma).\eeq
So, taking the derivative of the functional
$$(v,\varphi)\in H\times L^2(\Gamma;d\mu)\mapsto \int_\Gamma a(u(\gamma)-\varphi(\gamma) v,u(\gamma)-\varphi(\gamma) v;\gamma)d\mu(\gamma),$$
we deduce that if $(w,\psi)\in H\times L^2(\Gamma;d\mu)$ is a solution of (\ref{caruWbb}), with $w\not =0$, then
\beq\label{caruWcc}\psi(\gamma)={a(u(\gamma),w;\gamma)\over  a(w,w;\gamma)},\quad d\mu\hbox{-a.e. }\gamma\in \Gamma,\eeq
and  $w$ is a solution of the non-linear variational problem
\beq\label{caruWdd} \int_\Gamma{a(u(\gamma),w;\gamma)\over a(w,w;\gamma)}a(u(\gamma),v;\gamma)d\mu(\gamma)=\int_\Gamma{a(u(\gamma),w;\gamma)^2\over a(w,w;\gamma)^2}a(w,v;\gamma)d\mu(\gamma),\quad\forall\, v\in H.
\eeq
Note that if $w=0$, then $u=0$ and therefore $f=0$.\par
If $a$ does not depend on $\gamma$, statement (\ref{caruWdd}) can be written as
$$a\left(\int_\Gamma a(u(\gamma),w)\,u(\gamma)d\mu(\gamma),v\right)=a\left({\dis \int_\Gamma a(u(\gamma),w)^2d\mu(\gamma)\over a(w,w)}\, w,v\right),\quad\forall\, v\in H,$$
which implies that
$$\int_\Gamma a(u(\gamma),w)\,u(\gamma)d\mu(\gamma)={\dis \int_\Gamma a(u(\gamma),w)^2d\mu(\gamma)\over a(w,w)}\, w,$$
i.e. $w$ is an eigenvector of the operator
$$v\in H\mapsto \int_\Gamma a(u(\gamma),v)u(\gamma)d\mu(\gamma)$$
for the eigenvalue
$${\dis \int_\Gamma a(u(\gamma),w)^2d\mu(\gamma)\over a(w,w)}.$$
In contrast, when $a$ depends on $\gamma$ problem (\ref{caruWdd}) does not correspond to an eigenvalue equation.\end{remark}

  \clearpage

\section{One-dimensional approximations} \label{se:hgral}
In Section 4 we shall show the existence of the solution of problem (\ref{eq:pb1d}) for any arbitrary $k$. However a particularly interesting case from the point of view of the applications is $k=1$. We  dedicate this section to  this special case.
 Observe that for $Z \in {\cal S}_1$, there exists  $z \in H$ 
such that $Z=Span\{z\}$. 
The problem to solve can be  reformulated as follows. 
\begin{lemma} \label{eq:equivpbs}   Assume $f\not \equiv 0$. Then, the subspace $W \in \csu$ solves problem \eqref{eq:pb1dr} 
if and only if \, $W=\espan\{w\}$, where $w$ is a solution of
\beq\label{maxpb} 
 \max_{z\in H\atop z\not=0} \,\int_\Gamma \frac{\langle f(\gamma),z\rangle^2} {a(z,z;\gamma)}  \,d\mu(\gamma).
\eeq
\end{lemma}
\

\begin{proof} 
Let $Z\in \csu$. Then $Z=\espan \{z\}$, for some $z \in H$, and there exists a function 
$\varphi:\Gamma \mapsto \R$ such that
$$
\ugz=\varphi(\gamma)\, z,\quad d\mu\hbox{-a.e. }\gamma\in \Gamma.
$$
\par 
If $z\not =0$, then, as $\ugz$ is the solution to the variational equation (\ref{eq:pbgw2}),   we derive that
$$
\varphi(\gamma)=\displaystyle \frac{\langle f(\gamma),z \rangle}{a(z,z;\gamma)},\quad d\mu\hbox{-a.e. }\gamma\in \Gamma.
$$
 Using  this formula  we obtain that
\beq\label{foruno0}
\int_\Gamma\langle f,\ugz \rangle \,d \gamma=  \, \int_\Gamma \frac{\langle f(\gamma),z\rangle^2} {a(z,z;\gamma)}  \,d\mu(\gamma).
\eeq
\par If the maximum in (\ref{eq:pb1dr}) is obtained by a space of dimension one, then formula (\ref{foruno0}) proves  the desired result.\par
In contrast, if the maximum in (\ref{eq:pb1dr}) is obtained by the null space, then the maximum in ${\cal S}_1$ is equal to zero. Therefore the right-hand side of (\ref{foruno0}) is zero for every $z\in H$, which implies that $f=0$ $d\mu$-a.e. in $\Gamma$, in contradiction with the assumption $f\not \equiv 0$.
\end{proof}
\begin{remark}\label{rk:dir} Since the integrand which appears in (\ref{maxpb}) is homogenous of degree zero in $z$, problem \eqref{maxpb} is equivalent to
$$\max_{z\in H\atop \|z\|=1} \,\int_\Gamma \frac{\langle f(\gamma),z\rangle^2} {a(z,z;\gamma)}  \,d\mu(\gamma).
$$
\end{remark}
\
We now prove the  existence of a solution to problem \eqref{maxpb}. 
\begin{theorem} \label{th:main} Assume $f\not \equiv 0$.
Problem \eqref{maxpb} admits at least a solution.
\end{theorem}
Note that if $f \equiv 0$, then, every vector $w\in H\setminus\{0\}$ is a  solution of (\ref{maxpb}). \par\medskip\noindent
\begin{proof} Define
 \beq\label{max:M}
M^* :=  \sup_{z\in H\atop\|z\|=1} \,\int_\Gamma \frac{\langle f(\gamma),z\rangle^2} {a(z,z;\gamma)}  \,d\mu(\gamma),
\eeq
and consider a sequence $w_n\subset H$, with $\|w_n\|=1$ such that
\beq\label{conMas}
 \lim_{n\to\infty}\int_\Gamma \frac{\langle f(\gamma),w_n\rangle^2} {a(w_n,w_n;\gamma)}  \,d\mu(\gamma)=M^*.
\eeq
Up to a subsequence, we  can assume the existence of  $w\in H$, such that 
$w_n$  converges weakly in $H$ to  $w$. Taking into account that $f(\gamma)\in H'$, $a(\cdot,\cdot,\gamma)\in B_s(H)$ $d\mu$-a.e. $\gamma\in\Gamma$ and (\ref{eq:propab}) is satisfied, we get
\beq\label{fconv} \lim_{n\to\infty}\langle f(\gamma),w_n\rangle=\langle f(\gamma),w\rangle,\quad d\mu\hbox{-a.e. }\gamma\in\Gamma,\eeq
\beq\label{liminfa} \liminf_{n\to\infty}a(w_n,w_n;\gamma)\geq  a(w,w;\gamma),\quad d\mu\hbox{-a.e. }\gamma\in\Gamma.\eeq
On the other hand, we observe that (\ref{eq:propab}) and $\|w_n\|=1$ imply
\beq\label{aco1a}|\langle f(\gamma),w_n\rangle|\leq \|f(\gamma)\|_{H'},\ \ {1\over a(w_n,w_n;\gamma)}\leq {1\over\alpha}\qquad d\mu\hbox{-a.e. }\gamma\in\Gamma.\eeq
\par
If $w=0$, then (\ref{fconv}), (\ref{aco1a}) and Lebesgue's dominated convergence theorem imply
$$\lim_{n\to\infty}\int_\Gamma \frac{\langle f(\gamma),w_n\rangle^2} {a(w_n,w_n;\gamma)}  \,d\mu(\gamma)=0,$$
which by (\ref{conMas}) is equivalent to $M^\ast=0$. Taking into account  (\ref{eq:propab}) and the  definition (\ref{max:M}) of $M^\ast$, this is only possible if $f\equiv 0$ is the null function. As we are assuming $f\not \equiv 0$, we conclude that $ w$ is different of zero. Then, (\ref{aco1a}) proves
$$0\leq {\|f(\gamma)\|^2_{H'}\over \alpha}-{\langle f(\gamma),w_n\rangle^2\over a(w_n,w_n;\gamma)},\quad d\mu\hbox{-a.e. }\gamma\in\Gamma,$$
while (\ref{fconv}) and (\ref{liminfa})  prove
\begin{equation}\label{desFat}
\displaystyle\liminf_{n\to\infty} \left({\|f(\gamma)\|^2_{H'}\over \alpha}-{\langle f(\gamma),w_n\rangle^2\over a(w_n,w_n;\gamma)}\right)\geq {\|f(\gamma)\|^2_{H'}\over \alpha}-{\langle f(\gamma),w\rangle^2\over a(w,w;\gamma)},\quad d\mu\hbox{-a.e. }\gamma\in\Gamma.\eeq
Using   (\ref{conMas}), Fatou's lemma implies
$$\begin{array}{l}\displaystyle\int_\Gamma \left({\|f(\gamma)\|^2_{H'}\over \alpha}-{\langle f(\gamma),w\rangle^2\over a(w,w;\gamma)}\right)d\mu(\gamma)\leq \liminf_{n\to\infty}
\int_\Gamma \left({\|f(\gamma)\|^2_{H'}\over \alpha}-{\langle f(\gamma),w_n\rangle^2\over a(w_n,w_n;\gamma)}\right)d\mu(\gamma)\\ \noalign{\medskip}
\displaystyle=\int_\Gamma {\|f(\gamma)\|^2_{H'}\over \alpha}\,d\mu(\gamma)-M^\ast,
\end{array}
$$
or equivalently
\beq\label{desMas} M^\ast\leq \int_\Gamma {\langle f(\gamma),w\rangle^2\over a(w,w;\gamma)}d\mu(\gamma).\eeq
By definition (\ref{max:M}) of $M^\ast$, this proves that the above inequality is an equality and that $w$ is a solution of (\ref{maxpb}).
\end{proof}
\begin{remark} \label{metaproi} Actually, in place of (\ref{desFat}), one has the stronger result
$$\liminf_{n\to\infty} \left({\|f(\gamma)\|^2_{H'}\over \alpha}-{\langle f(\gamma),w_n\rangle^2\over a(w_n,w_n;\gamma)}\right)= {\|f(\gamma)\|^2_{H'}\over \alpha}-{\langle f(\gamma),w\rangle^2\over\displaystyle \liminf_{n\to\infty}a(w_n,w_n;\gamma))},\quad d\mu\hbox{-a.e. }\gamma\in\Gamma,$$
which by the proof used to prove 
 (\ref{desMas}) shows
 $$M^\ast\leq \int_\Gamma  {\langle f(\gamma),w\rangle^2\over \displaystyle\liminf_{n\to\infty}a(w_n,w_n;\gamma)}d\mu(\gamma).$$
Combined with 
$$M^\ast= \int_\Gamma {\langle f(\gamma),w\rangle^2\over a(w,w;\gamma)}d\mu(\gamma)$$
and (\ref{liminfa}), this implies
$$a(w,w;\gamma)=\liminf_{n\to\infty}a(w_n,w_n;\gamma)\quad d\mu\hbox{-a.e. }\gamma\in\Gamma\ \hbox{ such that }\langle f(\gamma),w\rangle\not =0.$$
By (\ref{eq:propab}) and $f\not \equiv 0$, this proves the existence of a subsequence of $w_n$ which converges strongly to $w$.\par
 Since this proof can be carried out by replacing $w_n$ by  any subsequence of $w_n$, we conclude that the whole sequence $w_n$ (which  we extracted just after (\ref{conMas}) assuming that it converges weakly to some $w$)
actually converges strongly to $w$.\par
The above result may be used to build a computable approximation of a solution of \eqref{maxpb}. Indeed, for $f\not \equiv 0$,
let $\{H_n\}_{n\ge 1}$ be an internal approximation of $H$, that is a sequence of subspaces  of finite dimension of $H$ such that 
$$
\dis\lim_{n\to \infty}  \inf_{\psi \in H_n} \|z - \psi\|=0,  \quad \forall z \in H.
$$ 
and consider a solution $w_n$ of
$$
 \max_{z\in H_n\atop \|z\|=1} \,\int_\Gamma \frac{\langle f(\gamma),z\rangle^2} {a(z,z;\gamma)}  \,d\mu(\gamma).
$$
The existence of such a $w_n$ can be obtained by reasoning as in the proof of Theorem
\ref{th:main} or just using Weierstrass theorem because  the dimension of $H_n$ is finite.\par
Taking $\tilde w$ a solution of (\ref{maxpb}) and a sequence $\tilde w_n\in H_n$ 
 converging to $\tilde w$ in $H$, we have
$$\begin{array}{l}\displaystyle
\int_\Gamma \frac{\langle f(\gamma),\tilde w\rangle^2} {a(\tilde w,\tilde w;\gamma)}  \,d\mu(\gamma)=\lim_{n\to\infty}\int_\Gamma \frac{\langle f(\gamma),\tilde w_n\rangle^2} {a(\tilde w_n,\tilde w_n;\gamma)}  \,d\mu(\gamma)\\ \noalign{\medskip}\displaystyle\leq\liminf_{n\to\infty}\int_\Gamma \frac{\langle f(\gamma),w_n\rangle^2} {a( w_n,w_n;\gamma)}  \,d\mu(\gamma)\leq \limsup_{n\to\infty}\int_\Gamma \frac{\langle f(\gamma),w_n\rangle^2} {a( w_n,w_n;\gamma)}  \,d\mu(\gamma)\leq \int_\Gamma \frac{\langle f(\gamma),\tilde w\rangle^2} {a(\tilde w,\tilde w;\gamma)}  \,d\mu(\gamma),\ea$$
and thus
$$ \lim_{n \to \infty} \int_\Gamma \frac{\langle f(\gamma),w_n\rangle^2} {a(w_n,w_n;\gamma)}  \,d\mu(\gamma)=\int_\Gamma \frac{\langle f(\gamma),\tilde w\rangle^2} {a(\tilde w,\tilde w;\gamma)}  \,d\mu(\gamma)=M^*.
$$
 This proves that the sequence $w_n$ satisfies (\ref{conMas}). Therefore any subsequence of $w_n$ which converges weakly to some $w$  converges strongly to $w$ which is a solution of \eqref{maxpb}.
\end{remark}

 \clearpage

\section{Higher-dimensional approximations} \label{se:hfdim}

This section is devoted to   the proof of the existence  of an optimal subspace which is 
solution of \eqref{eq:pb1d} when $k\geq 1$ is any given number.  

\begin{theorem} \label{ThexPk} For any given $k \ge 1$, problem  \eqref{eq:pb1d} admits at least one solution.
\end{theorem}

\begin{proof}
As in the proof of Theorem \ref{th:main},
we use the direct method of the Calculus of Variations. Denoting by $m_k$
\beq\label{defMk}m_k=\inf_{Z \in \csk}  \bar a(u-u_Z,u-u_Z),\eeq
we consider a sequence of spaces $W_n\in \csk$  such that  $w_n := u_{W_n}$ satisfies
\beq\label{conMk} \lim_{n\to\infty}\bar a(u-w_n,u-w_n)=m_k.\eeq
Taking into account that by Proposition \ref{procaruz}
\beq\label{inclZn} Z\subset \tilde Z\Longrightarrow \bar a(u-u_{\tilde Z},u-u_{\tilde Z})\leq \bar a(u-u_Z,u-u_Z),\eeq
we can assume that the dimension of $W_n$ is equal to $k$. Moreover, we observe  that  (\ref{conMk}) implies that $w_n$ is bounded in $L^2(\Gamma, H;d\mu)$.\par
 Let $(z_n^1,\cdots, z_n^k)$ be  an orthonormal basis of  $W_n$. It holds
\beq\label{descwn}
w_n(\gamma) = \sum_{j=1}^k (w_n(\gamma), z_n^j )\, z_n^j , \quad  d\mu\mbox{-a.e. } \gamma\in\Gamma.
\eeq
Since the norm of the vectors $z_n^j$ is one, there exists a subsequence of $n$ and $k$  vectors $z^j\in H$ such that
\beq\label{convpj} z_n^j\rightharpoonup z^j\ \hbox{ in }H,\quad\forall\,j\in\{1,\cdots,k\}.\eeq 
Using also
$$
|(w_n(\gamma), z_n^j )| \le  \|w_n(\gamma)\|, \qquad  d\mu\mbox{-a.e }\gamma\in \Gamma,
$$
we get that $(w_n,z_n^j)$ is bounded in $L^2(\Gamma,H;d\mu)$ for every $j$ and thus, there exists a subsequence of $n$ and $k$ functions  $p^j \in  L^2(\Gamma;d\mu)$ such that
\beq\label{conwnzn} (w_n,z_n^j)\rightharpoonup p^j\ \hbox{ in }L^2(\Gamma,H;d\mu),\ \forall\, j\in\{1,\cdots,k\}.\eeq 
\par
We claim that 
\beq\label{conwn}w_n\rightharpoonup w:=\sum_{j=1}^n p^jz^j \ \hbox{ in }L^2(\Gamma;d\mu).\eeq
Indeed, taking into account that  $w_n$ is bounded in $L^2(\Gamma, H;d\mu)$ and (\ref{descwn}), it is enough to show 
\beq\label{conwnznb}
\lim_{n \to \infty }  
\int_\Gamma \left((w_n, z_n^j )z_n^j , \varphi\,v \right)\, d\mu(\gamma)
 =\int_\Gamma (p^j z^j , \varphi\,v)\, d\mu(\gamma),\quad\forall\, \varphi\in L^2(\Gamma;d\mu),\ \forall\, v\in H.
\eeq
This is a simple consequence of 
$$
\int_\Gamma \left((w_n, z_n^j )z_n^j , \varphi\,v \right)\, d\mu(\gamma)=(z_n^j , v) 
 \int_\Gamma (w_n, z_n^j )\,\varphi\, d\mu(\gamma),$$
combined with (\ref{convpj}) and (\ref{conwnzn}).\par
From the continuity and convexity of the quadratic form associated to $\bar a$, as well as from (\ref{conwn}) and (\ref{conMk}), we  have
\beq\label{paliuwn}
\bar a(u-w,u-w)\le \lim_{n \to \infty} \bar a(u-w_n,u-w_n) = m_k.
\eeq
Using that $W={\rm Span}\{z^1,\cdots , z^k\} \in \mathcal S_k$, and that (see Proposition \ref{procaruz})
\beq\label{desiuW}
\bar a(u-u_W,u-u_W)\leq \bar a(u-w,u-w),
\eeq
we conclude that $W$ is a solution of \eqref{eq:pb1d}.
\end{proof}
\begin{remark} From (\ref{paliuwn}), (\ref{desiuW}), definition (\ref{defMk}) of $m_k$ and Proposition \ref{procaruz}, we have that $w=u_W$ in the proof of Theorem \ref{ThexPk}. Moreover, 
$$\bar a(u-w,u-w)=m_k=\lim_{n\to\infty}\bar a(u-w_n,u-w_n),$$
which combined with (\ref{conwn}) proves that $w_n$ converges strongly to $w$ in $L^2(\Gamma,H;d\mu)$. As in Remark \ref{metaproi}, this can be used to  build a strong approximation of a solution of \eqref{eq:pb1d} by using an internal approximation of $H$.
\end{remark}
\clearpage 
  \clearpage 
\section{An iterative algorithm by deflation}
In the previous section, for any given $k \ge 1$, we have proved the existence of an optimal subspace for problem (\ref{eq:pb1d}). We use here this fact to build an iterative approximation of the solution of (\ref{eq:pbg}) by a deflation approach. Let us denote 
\beq\label{dePiku}\Pi_k(v)=\left\{v_W\,|\quad
W\ \hbox{ solves }\ \min_{Z \in \csk}  \bar a(v-v_Z,v-v_Z)
\,\right\},\quad\forall\, v\in L^2(\Gamma,H;d\mu).\eeq
\par
The deflation algorithm is as follows

\begin{itemize}
\item Initialization: 
\begin{equation}\label{al:def1}
\dis u_0=0
\end{equation}
 \item Iteration: Assuming $u_{i-1} \in H$ known for $i =1,2,\cdots$, set 
\begin{equation}\label{al:def2}
u_i=u_{i-1} + s_i,\quad \mbox{ with }\, s_i \in \Pi_k(e_{i-1}), \,\mbox{ where   }e_{i-1}=u-u_{i-1}
\end{equation}
\end{itemize}
\begin{remark} Since one has $e_{i-1}=u-u_{i-1}$ by (\ref{al:def2}) and since by (\ref{dealu}), $u$ is the solution of 
(\ref{eq:pbg2}), the function $e_{i-1}$ satisfies
$$\left\{\ba{l}\dis e_{i-1}\in L^2(\Gamma,H;d\mu),\\ \ecart\dis \bar a(e_{i-1},v)=\int_\Gamma \langle f(\gamma),v(\gamma)\rangle\,d\mu(\gamma)-\bar a(u_{i-1},v),\ \ \forall\, v\in  L^2(\Gamma,H;d\mu).\ea\right.$$
Therefore Proposition \ref{pr:equiv} applied to the case where  $f$ is replaced by the function $\hat f_i$ defined by
$$\int_\Gamma\langle\hat f_i(\gamma),v(\gamma)\rangle d\mu(\gamma)=\int_\Gamma \langle f(\gamma),v(\gamma)\rangle\,d\mu(\gamma)-\bar a(u_{i-1},v),\ \ \forall\, v\in  L^2(\Gamma,H;d\mu),$$
proves that $s_i\in \Pi_k(e_{i-1})$  is equivalent to $s_i=(e_{i-1})_W$, where $W$ is a solution of
$$\max_{Z \in \csk}  \left\{\int_\Gamma \langle f(\gamma),(e_{i-1})_Z(\gamma)\rangle \, d\mu(\gamma)-\bar a(u_{i-1},(e_{i-1})_Z)\right\},$$
where, in accordance to (\ref{eq:pbgw2}), $(e_{i-1})_Z$ denotes the solution of 
$$\left\{\ba{l}\dis (e_{i-1})_Z\in L^2(\Gamma,Z;d\mu),\\ \ecart\dis
\bar a\big((e_{i-1})_Z,z\big)=\int_\Gamma\langle f(\gamma),z(\gamma)\rangle\, d\mu(\gamma)-\bar a(u_{i-1},z),\quad \forall\, z\in L^2(\Gamma,Z;d\mu).\ea\right.$$\par
This observation allows one to carry out the iterative process without knowing the function $u$ (compare with (\ref{al:def2})).
\end{remark}
\par
The convergence of the algorithm is given by the following theorem. Its proof follows the ideas of
\cite{falnou}.
\begin{theorem} \label{Thcdf} The sequence $u_i$ provided by the least-squares PGD algorithm (\ref{al:def1})-(\ref{al:def2}) strongly converges in $L^2(\Gamma,H;d\mu)$ to the parameterized solution $\gamma\in \Gamma \mapsto u(\gamma) \in H$ of problem (\ref{eq:pbg}).
\end{theorem}
\begin{proof} By (\ref{al:def2}) and Proposition \ref{pbeq1} applied to the case where $u$ is replaced by $e_{i-1}$, we have that $s_i$ is a solution of 
\beq\label{carsiW} \min_{v\in L^2(\Gamma,H;d\mu)\atop
{\rm dim}\,R(v)\leq k}\bar a(e_{i-1}-v,e_{i-1}-v).\eeq
This proves in particular that $s_i$ is a solution of
$$\min_{v\in L^2(\Gamma,H;d\mu)\atop R(v)\subset R(s_i)}\bar a(e_{i-1}-v,e_{i-1}-v),$$
and therefore
$$\bar a(e_{i-1}-s_i,v)=0,\quad\forall\, v\in L^2(\Gamma,H;d\mu)\ \hbox{ with }R(v)\subset R(s_i).$$
But (\ref{al:def2}) implies that
\beq\label{reeisi} e_{i-1}-s_i=e_i,\eeq
which gives
\beq\label{ortei}\bar a(e_i,v)=0,\quad\forall\, v\in L^2(\Gamma,H;d\mu)\ \hbox{ with }R(v)\subset R(s_i).\eeq
Taking  $v=s_i$  and using again (\ref{reeisi}) we get
\beq\label{aeieisisi}\bar a(e_{i-1},e_{i-1})=\bar a (s_i,s_i)+\bar a(e_i,e_i),\quad\forall\, i\geq 1,\eeq
and therefore
\beq\label{existim}\bar a(e_i,e_i)+\sum_{j=1}^i\bar a(s_j,s_j)=\bar a(e_0,e_0),\quad\forall\, i\geq 1.\eeq
Thus, we have
\beq\label{acsn}e_i\ \hbox{  is bounded in  }L^2(\Gamma,H;d\mu),\eeq
\beq\label{sumacot}\sum_{j=1}^\infty\bar a(s_j,s_j)\leq \bar a(e_0,e_0).\eeq\par
By (\ref{acsn}), there exists a subsequence  $e_{i_n}$ of $e_i$ and $e\in L^2(\Gamma,H;d\mu),$
such that
\beq\label{convsn}e_{i_n}\rightharpoonup e\ \hbox{ in }L^2(\Gamma,H;d\mu).\eeq
On the other hand, since $s_{i_n+1}$ is a solution of (\ref{carsiW}) with $i-1$ replaced by $i_n$, we get 
\beq\label{acotgm}\begin{array}{l}\displaystyle \bar a(e_{i_n}-s_{i_n+1},e_{i_n}-s_{i_n+1})\leq\bar  a(e_{i_n}-v,e_{i_n}-v)= \bar a(e_{i_n},e_ {i_n})-2\bar a(e_{i_n},v)+\bar a(v,v),\\ \ecart\dis
\forall\, v\in L^2(\Gamma,H;d\mu),\ {\rm dim}\,R(v)\leq k,\ea\eeq
and then
$$\bar a(e_{i_n}-s_{i_n+1},e_{i_n}-s_{i_n+1})- \bar a(e_{i_n},e_ {i_n})\leq -2\bar a(e_{i_n},v)+\bar a(v,v),\ \ \forall\, v\in L^2(\Gamma,H;d\mu),\ {\rm dim}\,R(v)\leq k,$$
or in other terms
$$-2\bar a(e_{i_n},s_{i_n+1})+\bar a(s_{i_n+1},s_{i_n+1})
\leq -2\bar a(e_{i_n},v)+\bar a(v,v),\ \ \forall\, v\in L^2(\Gamma,H;d\mu),\ {\rm dim}\,R(v)\leq k.$$
Thanks to (\ref{acsn}) and (\ref{sumacot}), the left-hand side tends to zero when $n$ tends to infinity, while in the right-hand side we can pass to the limit by (\ref{convsn}). Thus, we have
$$2\bar a(e,v)\leq \bar a(v,v),\quad\forall\, v\in L^2(\Gamma,H;d\mu),\ \ {\rm dim}\,{\rm R}(v)\leq k.$$
Replacing in this equality $v$ by $tv$ with $t>0$, dividing by $t$, letting $t$ tend to zero and writing the resulting inequality for $v$ and $-v$, we get
$$\bar a(e,v)=0,\quad\forall\, v\in L^2(\Gamma,H;d\mu),\ \ {\rm dim}\,{\rm R}(v)\leq k.$$\par
Taking $v=w\varphi$, with  $w\in H$, $\varphi\in L^2(\Gamma;d\mu)$, 
and recalling definition (\ref{defab}) of $\bar a$ we deduce
$$\int_\Gamma a(e(\gamma),w;\gamma)\,\varphi(\gamma)\,d\mu(\gamma)=0,\quad\forall\,z\in H,\ \forall\,\varphi\in L^2(\Gamma;d\mu),$$
and then  for any $w\in H$,  there exists a subset $N_w\in {\cal B}$ with $\mu(N_w)=0$ such that
$$a(e(\gamma),w;\gamma)=0,\quad \forall\,\gamma\in \Gamma\setminus N_w.$$
The separability of $H$ implies that $N_w$ can be chosen independent of $w$, and then we have
$$a(e(\gamma),w;\gamma)=0,\quad\forall\,w\in H,\ d\mu\hbox{-a.e. }\gamma\in\Gamma,$$
and therefore
\beq\label{idenli} e(\gamma)=0\quad d\mu\hbox{-a.e. }\gamma\in \Gamma.\eeq
This proves that $e$ does not depend on the subsequence in (\ref{convsn}) and that 
\beq\label{convsnb}e_i\rightharpoonup 0\ \hbox{ in }L^2(\Gamma,H;d\mu).\eeq\par
Let us now prove that in (\ref{convsnb}) the convergence is strong in $L^2(\Gamma,H;d\mu)$. We use that thanks to (\ref{reeisi}), we have
$$e_i=-\sum_{j=1}^i  s_j+e_0,\quad\forall\, i\geq 1,$$
and so,
\beq\label{acot1cf}\bar a(e_i,e_i)=-\sum_{j=1}^i\bar a(e_i,s_j)+\bar a(e_i,e_0),\quad\forall\,i\geq 1.
\eeq\par
In order to estimate the right-hand side of the latest equality, we introduce, for  $i,j\geq 1$, the function $z_{i,j}$ as the solution of
\beq\label{dezij} z_{i,j}\in  L^2(\Gamma,{\rm R}(s_j);d\mu),\qquad 
\bar a(z_{i,j},v)=\bar a(e_{i-1},v),\ \ \forall\, 
v\in  L^2(\Gamma,{\rm R}(s_j);d\mu).\eeq
We have
\beq\label{acoemzjm}\big|\bar a(e_{i-1},s_j)\big|=\big|\bar a(z_{i,j},s_j)\big|\leq \bar a(z_{i,j},z_{i,j})^{1\over 2}\,\bar a(s_j,s_j)^{1\over 2}.
\eeq
Using (\ref{aeieisisi}), (\ref{reeisi}), the fact that $s_i $ is a solution of (\ref{carsiW}) and dim$\,{\rm R}(s_j)\leq k$ 
$$\bar a(e_{i-1},e_{i-1})-\bar a(s_i,s_i)=\bar a(e_{i-1}-s_i,e_{i-1}-s_i)\leq\bar a(e_{i-1}-z_{i,j},e_{i-1}-z_{i,j}).$$
Expending the right-hand side and using $v=z_{i,j}$ in (\ref{dezij}) this gives
$$\bar a(z_{i,j},z_{i,j})\leq \bar a(s_i,s_i),$$
which combined with (\ref{acoemzjm}) provides the estimate 
$$\big|\bar a(e_{i-1},s_j)\big|\leq \bar a(s_i,s_i)^{1\over 2}\,\bar a(s_j,s_j)^{1\over 2},\quad\forall\, i,j\geq 1.$$
\par
Using the latest estimate in  (\ref{acot1cf}) and then Cauchy-Schwarz's inequality, we get
\beq\label{estfin}\left\{\ba{l}\dis\bar a(e_i,e_i)\leq  \bar a(s_{i+1},s_{i+1})^{1\over 2}\sum_{j=1}^i\bar a(s_j,s_j)^{1\over 2}+\bar a(e_i,e_0)\\ \ecart\dis
\leq  \bar a(s_{i+1},s_{i+1})^{1\over 2}\,i^{1\over 2}\left(\sum_{j=1}^\infty\bar a(s_j,s_j)\right)^{1\over 2}+\bar a(e_i,e_0),\quad\forall\,i\geq 1.\ea\right.\eeq
But the criterion of comparison of two series with nonnegative terms and the facts that (see (\ref{sumacot}))
$$\sum_{i=1}^\infty{1\over i}=\infty,\qquad \sum_{i=1}^\infty\bar a(s_i,s_i)<\infty,$$ prove that
$$\liminf_{i\to\infty}\,\bar a(s_{i+1},s_{i+1})\,i=\liminf_{i\to\infty}{\bar a(s_{i+1},s_{i+1})\over {1\over i}}=0.$$
Since $\bar a(e_i,e_i)$ is a decreasing sequence by (\ref{aeieisisi}) and since (\ref{convsnb})  asserts that $e_i$ converges weakly to zero, we can pass to the limit in  (\ref{estfin}), to deduce
$$\ba{l}\displaystyle\lim_{i\to\infty}\bar a(e_i,e_i)=\liminf_{i\to\infty}\,\bar a(e_i,e_i)\\ \ecart\dis\leq  \liminf_{i\to\infty}\left(\,  \bar a(s_{i+1},s_{i+1})^{1\over 2}\,i^{1\over 2}\left(\sum_{j=1}^\infty\bar a(s_j,s_j)\right)^{1\over 2}+\bar a(e_i,e_0)\right)=0.\ea$$
This proves that $e_i$ converges strongly to zero in $L^2(\Gamma,H;d\mu)$. Since $e_i=u-u_i$ this finishes the proof of Theorem \ref{Thcdf}. 
\end{proof}
\begin{remark}
In many cases the corrections $s_i$ decrease exponentially in the sense that: 
$$
\|s_i\| =O(\rho^{-i})\,\, \mbox{ as  } i \to +\infty,\ \hbox{ for some }\rho>1.
$$
This occurs in particular for the standard POD expansion when $\Gamma$ is an open set of $\R^N$, $\mu$ is the Lebesgue measure and the function $f=f(\gamma)$ is analytic with respect to $\gamma$ (see \cite{azbelcha}).  Then $\|s_{i}\| $ is a good estimator for the error $\|u-u_i\|$.
\end{remark}
  \clearpage 
\section{Conclusion}
In this paper we have introduced an iterative deflation algorithm to solve parametric symmetric elliptic equations. It is a Proper Generalized Decomposition algorithm as it builds a tensorized representation of the parameterized solutions, by means of optimal subspaces that minimize the residual in mean quadratic norm. It is intrinsic in the sense that in each deflation step the residual is minimized in the \lq\lq natural" mean quadratic norm generated by the parametric elliptic operator. It is conceptually close to the Proper Orthogonal Decomposition with the difference that in the POD the residual is minimized with respect to a fixed mean quadratic norm. Due to this difference, spectral theory  cannot be applied.

We have proved the existence of the optimal subspaces of dimension less than or equal to a fixed number, as required in each iteration of the deflation algorithm, with a specific analysis for the one-dimensional case. Also, we have proved the strong convergence in the natural mean quadratic norm of the deflation algorithm for quite general parametric elliptic operators. 

We will next focus our research on the analysis of convergence rates of the deflation algorithm that we introduced. We will compare the convergence rates with those of the POD expansion, to determine whether the use of the \lq\lq natural" mean quadratic norm provides improved convergence rates. We will also work on the numerical approximation of the algorithm, based upon a \lq\lq trust" solution on high-fidelity finite-dimensional subspaces of the given Hilbert space, to construct a feasible Reduced Order Modeling algorithm.\par
All the results obtained in the present paper refer to $a$ symmetric. In a future work we will consider the non-symmetric case.

\section*{Acknowledgements}

The work of Mejdi Azaiez and Tom\'as Chac\'on has been partially supported by the Spanish Government - Feder EU grant MTM2015-64577-C2-1-R.\par\medskip
The work of Juan Casado-D\'{\i}az and Fran\c cois Murat  has been partially supported by the Spanish Government - Feder EU grant MTM2014-53309-P.

  \clearpage



\end{document}